\input epsf 
\documentclass[a4paper,11pt]{amsart}
\title{Hooks and powers of parts in partitions}
\author{Roland Bacher and Laurent Manivel}
\date{last modified Aug 29, 2001}

\begin{document}
\maketitle

\vskip0.5cm
{\it Abstract\footnote {Math. Class.: 05A17 Keywords: Partition.}: 
This paper shows that the number of hooks
of length $k$ contained in all partitions of $n$ equals $k$ times the 
number of parts of length $k$ in partitions of $n$. It contains also 
formulas for the moments (under uniform distribution) of
$k-$th parts in partitions of $n$.}

\section{Introduction and Main Results}     

Many textbooks contain material on partitions. Two standard references are
{\bf [A]} and {\bf  [S]}.

A {\it partition} of a natural integer $n$ with {\it parts}
$\lambda_1,\dots,\lambda_k$ is a finite decreasing
sequence $\lambda=(\lambda_1\geq \lambda_2\geq \dots\geq \lambda_k>0)$
of natural integers $\lambda_1,\dots ,\lambda_k>0$ such that $n=\sum_{i=1}^k
\lambda_i$. We denote by $\vert \lambda \vert$ the 
{\it content} $n$ of $\lambda$.
Partitions are also written as sums: $n=\lambda_1+\dots+\lambda_k$
and one uses also the (abusive) multiplicative notation
$$\lambda=1^{\nu_1}\cdot 2^{\nu_2}\cdots n^{\nu_n}$$
where $\nu_i$ denotes the number of parts equal to $i$ in the
partition $\lambda$. 

A partition is graphically represented by its {\it Young diagram} obtained
by drawing $\lambda_1$ adjacent boxes of identical size
on a first row, followed by $\lambda_2$ adjacent boxes of identical size
on a second row and so on with all first boxes (of different rows)
aligned along a common first column. In the sequel we identify a 
partition with its Young diagram. A {\it hook} in a partition is a choice
of a box $H$ in the corresponding Young diagram together with all boxes
at the right of the same row and all boxes below of the same column.
The total number of boxes in a hook is its {\it hooklength}, the number
of boxes in a hook to the right of $H$ is its {\it armlength}
and the number of boxes of a hook below $H$ is the {\it leglength}.
The Figure below displays the Young diagram of the partition $(5,4,3,1)$
of $13$ together with a hook of length $4$ having armlength $2$ and
leglength $1$. We call the couple $(\hbox{armlength,leglength})$
of a hook its {\it hooktype} and denote it by $\tau=\tau(\alpha,k-1-\alpha)$
if its armlength is $\alpha$ and its leglength $k-1-\alpha$. Such a hook has
hence total length $k$ and there are exactly $k$ different hooktypes
for hooks of length $k$. 

\medskip
\centerline{\epsfysize2cm\epsfbox{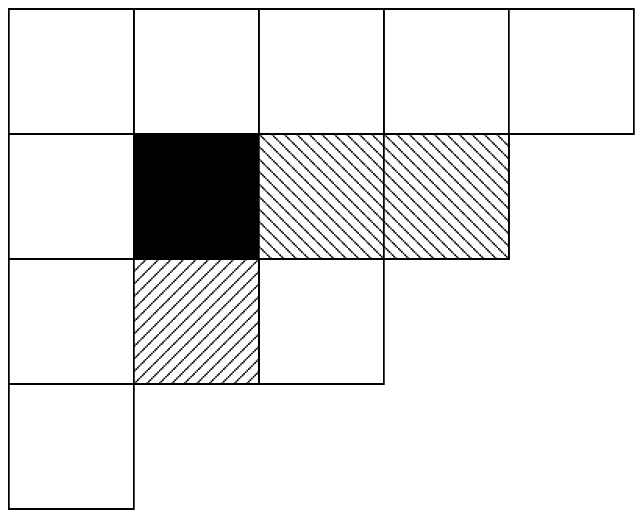}}
\centerline{The partition $(5,4,3,1)$ of $13$ together with a 
hook of type $\tau(2,1)$ and length $4$.}
\bigskip
Let $k$ be a natural integer and let $\tau=\tau(\alpha,k-1-\alpha)$ be the
hooktype of a hook of length $k$
with armlength $\alpha$ and leglength $k-1-\alpha$.
Given a partition $\lambda$ of $n$, set
$$\tau(\lambda)=\sharp\{\hbox{hooks of type }\tau\hbox{ in 
(the Young diagram of) }\lambda\}$$
and
$$\tau(n)=\sum_{\lambda,\ \vert \lambda\vert=n} \tau(\lambda)$$
where the sum is over all partitions of $n$.

{\bf Theorem 1.1.} {\sl One has
$$\begin{array}{lcl}
\displaystyle \sum_{n=1}^\infty\tau(n)z^n&\displaystyle =&\displaystyle
\frac{z^k}{1-z^k}\prod_{i=1}^\infty\frac{1}{
1-z^i}\cr&\displaystyle
=&\displaystyle
\sum_{\lambda=1^{\nu_1}2^{\nu_2}\cdots} \nu_k\ z^{\vert\lambda\vert}
\end{array}$$
where the last sum is over all partitions of integers.}

In other terms, the number of hooks of given type and length $k$
appearing in all partitions of $n$ equals the number of parts of length $k$ in
all partitions of $n$.

This result implies in particular that the total number of hooks of given
type $\tau=\tau(\alpha,k-1-\alpha)$
occuring in all partitions of $n$ depends only on the length $k$ 
and not on the particular hooktype $\tau(\alpha,k-1-\alpha)$ itself. Since there
are exactly $k$ distinct hooktypes for hooks of length $k$, the total
number of hooks of length $k$ in partitions of $n$ is given by
the coefficient of $z^n$ of the series 
$$k\frac{z^k}{1-z^k}\prod_{i=1}^\infty\frac{1}{1-z^i}\ .$$

For $\alpha,\beta\geq 0$ two natural integers, define the {\it $q-$binomial
${\alpha+\beta\choose \alpha}_q$} by
$${\alpha+\beta\choose \alpha}_q=\frac{\prod_{j=1}^{\alpha+\beta}(q^j-1)}{
\prod_{j=1}^{\alpha}(q^j-1)\prod_{j=1}^{\beta}(q^j-1)}\ .$$
The $q-$binomial coefficient ${\alpha+\beta\choose \alpha}_q$ is a
polynomial of degree $\alpha\beta$ in $q$ with the coefficient of $q^n$
enumerating all partitions of $n$ having at most $\beta$ non-zero parts
which are all of length at most $\alpha$. 

The main ingredient of the proof of Theorem 1.1 is the following result
which is perhaps also of independent interest.

{\bf Proposition 1.2.} {\sl One has for $\alpha,\ \beta\in {\bf N}$
$${\alpha+\beta\choose \alpha}_q=\frac{1}{1-q^{\alpha+\beta+1}}\left(\sum_{i=0}^\infty
q^{i(\beta+1)}\prod_{j=i+1}^{i+\alpha}(1-q^j)\right)^{-1}\ .$$}

As remarked previously, partitions of a natural integer $n$ can 
be written in (at least) two different ways: either by considering
the finite decreasing sequence
$$\lambda=(\lambda_1\geq \lambda_2\geq \dots \geq \lambda_k)$$
of its parts or by considering the vector 
$$\nu=(\nu_1,\nu_2,\dots,\nu_n)$$
where $\nu_i$ counts the multiplicity of parts with length $i$ in $\lambda$.
We still denote by $\lambda$ the vector 
$$(\lambda_1,\lambda_2,\dots\lambda_k,0,\dots,0)\in {\bf Z}^n$$
of length $n$ obtained by appending $n-k$ zero-coordinates at the end of 
the vector $(\lambda_1,\dots,\lambda_k)$ defining a partition $\lambda$.

We consider moreover the vector $\gamma=(\gamma_1,\dots,\gamma_n)$
where $\gamma_i$ equals the number of coordinates equal to $i$ among
$\nu_1,\dots,\nu_n$. The vector $\gamma$ of a partition encodes
\lq\lq multiplicities of multiplicities (of parts)'' and does no longer
encode the partition since for instance
the $\gamma-$vectors of the two partitions
$5=3+2$ and $5=4+1$ give both rise to $\gamma=(2,0,0,0,0)$. 

We introduce also the vectors 
$$\begin{array}{l}
\lambda(n)=(\lambda_1(n),\lambda_2(n),\dots,\lambda_n(n))=\sum_{\vert\lambda\vert=n}
\lambda\ ,\cr
\nu(n)=(\nu_1(n),\nu_2(n),\dots,\nu_n(n))=\sum_{\vert\lambda\vert=n}
\nu\ ,\cr
\gamma(n)=(\gamma_1(n),\gamma_2(n),\dots,\gamma_n(n))=\sum_{\vert\lambda\vert=n}
\gamma\end{array}$$
of ${\bf Z}^n$ obtained by summing up the vectors $\lambda,\nu$  or $\gamma$
over all partitions of $n$. The coordinates of the vector $\lambda(n)$ are
of course related to the mean length (under uniform distribution) of the 
$k-$th part in partitions of $n$. Similarly, coordinates of $\nu(n)$ 
relate to the mean multiplicity of parts equal to $k$ and coordinates of 
$\gamma(n)$ measure the mean number of distinct part-lengths appearing with 
common multiplicity $k$.

The following Table displays all five partitions of $4$ together with
the corresponding $\lambda-,\ \nu-$ and $\gamma-$vectors.

$$\begin{array}{lccccc}
\hbox{Partition}&(1,1,1,1)&(2,1,1)&(2,2)&(3,1)&(4)\cr
\lambda-\hbox{vector}&(1,1,1,1)&(2,1,1,0)&(2,2,0,0)&(3,1,0,0)&(4,0,0,0)\cr
\nu-\hbox{vector}&(4,0,0,0)&(2,1,0,0)&(0,2,0,0)&(1,0,1,0)&(0,0,0,1)\cr
\gamma-\hbox{vector}&(0,0,0,1)&(1,1,0,0)&(0,1,0,0)&(2,0,0,0)&(1,0,0,0)
\end{array}$$

Summing up all $\lambda-,\ \nu-$ and $\gamma-$vectors associated to
the five partitions of $4$ we get hence
$$\begin{array}{l}
\lambda(4)=(12,5,2,1),\cr
\nu(4)=(7,3,1,1),\cr
\gamma(4)=(4,2,0,1)\ .\end{array}$$

One has the following result.

{\bf Theorem 1.3.} {\sl For all $n\geq 1$ we have $\lambda_n(n)=
\nu_n(n)=\gamma_n(n)=1$ and
$$\begin{array}{l}
\lambda_k(n)=\nu_k(n)+\lambda_{k+1}(n)\ ,\cr
\nu_k(n)=\gamma_k(n)+\nu_{k+1}(n)\end{array}$$
for $k=1,\dots,n-1.$}

These equalities can be restated as
$$\lambda_k(n)=\sum_{i=k}^n\nu_i(n)\hbox{ and }\nu_k(n)=\sum_{i=k}^n \gamma_k(n)\ .$$
This last equality states for instance that the sum over all partitions of $n$ 
of the number
of distinct parts arising with multiplicity at least $k$ equals the number of parts
equal to $k$ in all partitions of $n$. Our proof of this fact uses generating 
series. It would be interesting to have a direct combinatorial proof
of this.

The coordinates of the vectors $\nu(n)$ are of course given by the 
generating series mentionned in Theorem 1.1, ie. the $k-$th 
coordinate $\nu_k(n)$ of $\nu(n)$ equals the coefficient of $z^n$
in the generating series
$$\frac{z^k}{1-z^k}\prod_{i=1}\frac{1}{1-z^i}\ .$$

The coordinates of the vectors $\lambda(n)$ and $\nu(n)$ are then easily computed
using Theorem 1.3. More precisely, one has the following result:

{\bf Corollary 1.4.} {\sl (i) The $k-$th coordinate of the vector 
$\lambda(n)$ is the coefficient of $z^n$ in the generating series
$$\prod_{i=1}^\infty\frac{1}{1-z^i}\ \sum_{j=k}^\infty 
\frac{z^j}{1-z^j}\ .$$

\ \ (ii) The $k-$th coordinate of the vector 
$\gamma(n)$ is the coefficient of $z^n$ in the generating series
$$\frac{(1-z)\ z^k}{(1-z^k)(1-z^{k+1})}\prod_{i=1}^\infty\frac{1}{1-z^i}\ .$$
}

Given a partition 
$$\lambda=(\lambda_1,\dots,\lambda_n)=(1^{\nu_1}\cdots n^{\nu_n})$$
and an integer $d\geq 0$ we introduce the vectors ${\lambda\choose d}$
and ${\nu\choose d}\in{\bf Z}^n$ by setting
$${\lambda\choose d}=\left({\lambda_1\choose d},\dots,{\lambda_n\choose d}\right)
\hbox{ and }{\nu\choose d}=\left({\nu_1\choose d},\dots,{\nu_n\choose d}\right)$$
and define
$${\lambda(n)\choose d}=\sum_{\vert\lambda\vert=n}{\lambda\choose d}\ ,\quad
{\nu(n)\choose d}=\sum_{\vert 1^{\nu_1}\ 2^{\nu_2}\cdots\vert=n}{\nu\choose d}$$
with coordinates
$${\lambda_k(n)\choose d}=\sum_{\vert\lambda\vert=n}{\lambda_k\choose d}\ ,\quad
{\nu_k(n)\choose d}=\sum_{\vert 1^{\nu_1}\ 2^{\nu_2}\cdots\vert=n}{\nu_k\choose d}
\ .$$

The following example shows the vectors  ${\lambda\choose 1}=\lambda,
{\lambda\choose 2},{\lambda\choose 3}$ and ${\nu\choose 1}=\nu,{\nu\choose 2},
{\nu\choose 3}$ associated to all five partitions of $4$.

{\bf Example:}
$$\begin{array}{lcccccc}
\hbox{Partition}&&(1,1,1,1)&(2,1,1)&(2,2)&(3,1)&(4)\cr
{\lambda\choose 1}& =&(1,1,1,1)&(2,1,1,0)&(2,2,0,0)&(3,1,0,0)&(4,0,0,0)\cr
{\lambda\choose 2}& =&(0,0,0,0)&(1,0,0,0)&(1,1,0,0)&(3,0,0,0)&(6,0,0,0)\cr
{\lambda\choose 3}& =&(0,0,0,0)&(0,0,0,0)&(0,0,0,0)&(1,0,0,0)&(4,0,0,0)\cr
\end{array}$$
$$\begin{array}{lcccccc}
\hbox{Partition}&&(1,1,1,1)&(2,1,1)&(2,2)&(3,1)&(4)\cr
{\nu\choose 1}& =&(4,0,0,0)&(2,1,0,0)&(0,2,0,0)&(1,0,1,0)&(0,0,0,1)\cr
{\nu\choose 2}& =&(6,0,0,0)&(1,0,0,0)&(0,1,0,0)&(0,0,0,0)&(0,0,0,0)\cr
{\nu\choose 3}& =&(4,0,0,0)&(0,0,0,0)&(0,0,0,0)&(0,0,0,0)&(0,0,0,0)\cr
\end{array}$$

We have thus
$$\begin{array}{lll}
{\lambda(4)\choose 1}=(12,5,2,1),\quad&
{\lambda(4)\choose 2}=(11,1,0,0),\quad& 
{\lambda(4)\choose 3}=(5,0,0,0)\cr
{\nu(4)\choose 1}=(7,3,1,1),\quad&
{\nu(4)\choose 2}=(7,1,0,0),\quad& 
{\nu(4)\choose 3}=(4,0,0,0)\end{array}$$

The following probably well-known result allows easy computations of the
vectors ${\lambda(n)\choose d}$ and ${\nu(n)\choose d}$.

{\bf Proposition 1.5.} {\it For any natural integer $d\geq 0$, the
$k-$th coefficients ${\lambda_k(n)\choose d}$, respectively
${\nu_k(n)\choose d}$ (extended by ${0\choose d}$ for $k>n$) have 
generating series
$$\sum_{n}^\infty {\lambda_k(n)\choose d}\ z^n=
\left(\prod_{j=1}^{k-1}\frac{1}{1-z^j}\right)\left(\sum_{i=0}^\infty
{i\choose d}\ z^{ik}\ \left(\prod_{j=1}^i\frac{1}{1-z^j}\right)\right)$$
and 
$$\sum_{n}^\infty {\nu_k(n)\choose d}\ z^n=\left(\frac{z^k}{1-z^k}\right)^d
\left(\prod_{j=1}^{\infty}\frac{1}{1-z^j}\right)\ .$$
}

{\bf Remark 1.6.} One has
$$\sum_{\vert\lambda\vert=n} \lambda_k^d=\sum_i i!\ \hbox{Stirling}_2(d,i)\
{\lambda_k(n)\choose i}$$
and
$$\sum_{\vert 1^{\nu_1}\ 2^{\nu_2}\cdots\vert=n} \nu_k^d=
\sum_i i!\ \hbox{Stirling}_2(d,i)\
{\nu_k(n)\choose i}$$
where $\hbox{Stirling}_2(d,i)$ denote Stirling numbers of the second kind, 
defined by
$x^d=\sum_i\hbox{Stirling}_2(d,i)\ x(x-1)\cdots(x-i+1)$.

Asymptotics are not so easy to work out from 
the formula for ${\lambda(n)\choose d}$. 
Our last result is an equivalent expression
for the above series on which asymptotics are easier to see.

We introduce the generating series $\sigma_r(k)$ defined as
$$\sigma_r(k)=\sum_{i=k}^\infty \left(\frac{z^i}{1-z^i}\right)^r$$
for $r\geq 1$ and $k\geq 1$ natural integers. We consider the series
$\sigma_r(k)$ as beeing graded of degree $r$ and define the homogeneous series
$S_d(k)$ of degree $d$ by
$$\begin{array}{rcl}
\displaystyle S_d(k)&\displaystyle =&
\displaystyle \sum_{\vert(1^{\nu_1}\ 2^{\nu_2}\cdots)\vert=d}
\frac{d!}{\left(\sum_i \nu_i\right)!}{\left(\sum_i \nu_i\right)
\choose \nu_1\ \nu_2\dots}\prod_{i=1}^d \left(\frac{\sigma_i(k)}{i}
\right)^{\nu_i}\cr
&\displaystyle =&\displaystyle d!\quad\sum_{\vert 1^{\nu_1}2^{\nu_2}\cdots
t^{\nu_t}\vert=d}\qquad
\prod_{j=1}^d \frac{(\sigma_j(k))^{\nu_j}}{j^{\nu_j}\ \nu_j!}\end{array}$$
(ie. the coefficient of the homogeneous \lq\lq monomial'' series
$\sigma_\lambda(k)=\sigma_{\lambda_1}(k)
\dots \sigma_{\lambda_s}(k)$ equals the number of elements in the symmetric group
on $\vert \lambda \vert$ elements of the conjugacy class with cycle structure
$\lambda=(\lambda_1,\lambda_2,\dots,\lambda_s)$).

We have then the following result.

{\bf Theorem 1.7.} {\it For any natural integers $d\geq 1$ and $k\geq 1$,
we have
$$\sum_{n}^\infty {\lambda_k(n)\choose d}\ z^n=
\frac{S_d(k)}{d!}\ \left(\prod_{j=1}^{\infty}\frac{1}{1-z^j}\right)\ .$$}


The first series $S_i=S_i(k)$ are given in terms of $\sigma_j=\sigma_j(k)$ as follows
$$\begin{array}{l}
S_0=1,\cr
S_1=\sigma_1,\cr
S_2=\sigma_1^2+\sigma_2,\cr
S_3=\sigma_1^3+3\sigma_1\sigma_2+2\sigma_3\cr
S_4=\sigma_1^4+6\sigma_1^2\sigma_2+3\sigma_2^2+8\sigma_1\sigma_3+6\sigma_4\cr
S_5=\sigma_1^5+10\sigma_1^3\sigma_2+20\sigma_1^2\sigma_3+15\sigma_1\sigma_2^2+
30\sigma_1\sigma_4+20\sigma_2\sigma_3+24\sigma_5\cr\end{array}$$

Let us remark that the analogous statement of the Theorem 1.9 
for the generating series $\sum_n {\nu_k\choose d} z^n$ boils down to 
a trivial identity.

The formulas of Theorem 1.7 ease the computations of asymptotics (in $n$) for 
$\lambda_k(n)$ and its moments and allow probably a rederivation of the results
contained in {\bf [EL]} and {\bf [VK]}

\section{Proofs}

{\bf Proof of Proposition 1.2.} 
Since ${\beta\choose 0}_q=1$, for $\alpha=0$ the proposition
boils down to the well-known formula for the geometric series
$$\sum_{i=0}^\infty q^{i(\beta+1)}=\frac{1}{1-q^{\beta+1}}\ .$$

The proposition is equivalent to the identity
$$\left(\prod_{j=\alpha+1}^{\alpha+\beta+1} (1-q^j)\right)\left(
\sum_{i=0}q^{i(\beta+1)}\prod_{j=i+1}^{i+\alpha}(1-q^j)\right)=\prod_{j=1}^\beta
(1-q^j)\ .$$
Since the right-hand side of this expresssion depends only on $\beta$,
it is enough to show that the expression
$$\left(\prod_{j=\alpha+1}^{\alpha+\beta+1}(1-q^j)\right)\left(\sum_{i=0}q^{i(\beta+1)}
\prod_{j=i+1}^{i+\alpha}(1-q^j)\right)$$
$$-\left(\prod_{j=\alpha}^{\alpha+\beta}(1-q^j)\right)
\left(\sum_{i=0}q^{i(\beta+1)}
\prod_{j=i+1}^{i+\alpha-1}(1-q^j)\right)$$
equals zero for all $\alpha\geq 1$.

Dividing by $\prod_{j=\alpha+1}^{\alpha+\beta}(1-q^j)$ we 
get
$$(1-q^{\alpha+\beta+1})\left(\sum_{i=0}q^{i(\beta+1)}
\prod_{j=i+1}^{i+\alpha}(1-q^j)\right)-
(1-q^\alpha)\left(\sum_{i=0}q^{i(\beta+1)}
\prod_{j=i+1}^{i+\alpha-1}(1-q^j)\right)$$
$$=\sum_{i=0}q^{i(\beta+1)}
\prod_{j=i+1}^{i+\alpha-1}(1-q^j)-\sum_{i=0}q^{i+\alpha+i(\beta+1)}
\prod_{j=i+1}^{i+\alpha-1}(1-q^j)$$
$$-q^{\alpha+\beta+1}\sum_{i=0}q^{i(\beta+1)}(1-q^{i+\alpha})
\prod_{j=i+1}^{i+\alpha-1}(1-q^j)$$
$$-\sum_{i=0}q^{i(\beta+1)}\prod_{j=i+1}^{i+\alpha-1}(1-q^j)+
q^\alpha\sum_{i=0}q^{i(\beta+1)}\prod_{j=i+1}^{i+\alpha-1}(1-q^j)$$
$$=q^\alpha\left(\sum_{i=0}(1-q^i)q^{i(\beta+1)}
\prod_{j=i+1}^{i+\alpha-1}(1-q^j)-q^{\beta+1}(1-q^{i+\alpha})q^{i(\beta+1)}
\prod_{j=i+1}^{i+\alpha-1}(1-q^j)\right)$$
$$=q^\alpha\left(\sum_{i=0}q^{i(\beta+1)}
\prod_{j=i}^{i+\alpha-1}(1-q^j)-\sum_{i=0}q^{(i+1)(\beta+1)}
\prod_{j=(i+1)}^{(i+1)+\alpha-1}(1-q^j)\right)=0$$
which proves Proposition 1.2.\hfill QED

{\bf Proof of Theorem 1.1.} 
Let us consider a hooktype $\tau=\tau(\alpha,k-1-\alpha)$ 
of length $k$ with armlength $\alpha$ and leglength $k-1-\alpha$.
The number 
$$\tau(n)=\sum_{\vert\lambda\vert=n} \tau(\lambda)$$
equals the coefficient of $z^n$ in the generating series
$$\sum_{i=0}\left(\prod_{j=1}^i\frac{1}{1-z^j}\right)
z^{(i+1)(k-\alpha)+\alpha}P_{k-1-\alpha,\alpha}(z)\prod_{j=i+1+\alpha}^\infty
\frac{1}{1-z^j}$$
$$=\left(\prod_{j=1}^\infty \frac{1}{1-z^j}  \right) P_{k-1-\alpha,\alpha}(z) 
\ z^k
\left(\sum_{i=0}z^{i(k-\alpha)}\prod_{j=i+1}^{i+\alpha}(1-z^j)\right)$$
(the factor $\left(\prod_{j=1}^i\frac{1}{1-z^j}\right)$ accounts for all 
what happens below a hook $H$ of type $\tau$, the factor 
$z^{(i+1)(k-\alpha)+\alpha}P_{k-1-\alpha,\alpha}(z)$ accounts for rows involved
in $H$ and the last factor $\prod_{j=i+1+\alpha}^\infty
\frac{1}{1-z^j}$ depends on what happens on rows above $H$)
where $P_{\alpha,\beta}(z)$ denotes the generating series of partitions
having at most $\alpha$ parts and all parts are of length at most $\beta$.
The generating series $P_{\alpha,\beta}(q)$ is by definition
the $q-$binomial coefficient
$$P_{\alpha,\beta}(q)={\alpha+\beta\choose \alpha}_q=
\frac{\prod_{j=1}^{\alpha+\beta}(q^j-1)}
{\prod_{j=1}^{\alpha}(q^j-1)\prod_{j=1}^{\beta}(q^j-1)}\ .$$

Applying Proposition 1.2 with $\beta=k-1-\alpha$ we have
$$\sum_n \tau(n)z^n=
\left(\prod_{j=1}^\infty\frac{1}{1-z^j}\right)\cdot $$
$$\cdot \frac{1}{1-z^k}\left(\sum_{i=0}
z^{i(k-\alpha)}\prod_{j=i+1}^{i+\alpha}(1-z^j)\right)^{-1}z^k
\left(\sum_{i=0}z^{i(k-\alpha)}\prod_{j=i+1}^{i+\alpha}(1-z^j)\right)$$
$$=\frac{z^k}{1-z^k}\prod_{j=1}^\infty\frac{1}{1-z^j}$$
which proves the first equality of the Theorem. 

The last equality follows from the easy identities
$$\sum_{\lambda=(1^{\nu_1}\ 2^{\nu_2}\cdots)} \nu_k\ z^{\vert\lambda\vert}=
\sum_{i\geq 1} iz^{ik}\ \prod_{1\leq j\not= i} \frac{1}{1-z^j}=
\frac{z^k}{1-z^k}\prod_{j= 1}^\infty \frac{1}{1-z^j}$$
thus finishing the proof.\hfill QED

{\bf Proof of Theorem 1.3.} The partition $1^n$ yields the unique non-zero 
contribution to $\lambda_n(n)$ and $\gamma_n(n)$ and this contribution equals 
one in both cases. The partition $n$ consisting of a unique part of
length $n$ yields the unique non-zero contribution to $\nu_n(n)$ and
this contribution equals again $1$.

Given a partition 
$\lambda=(\lambda_1,\dots,\lambda_k)$ of $n$, the {\it conjugate
partition} $\lambda^t=(\lambda^t_1,\dots,\lambda^t_{k'})$ of $\lambda$
is defined by
$$\lambda^t_j=\sharp\{i\ \vert\ \lambda_i\geq j\}$$
(this corresponds to a reflection of the Young diagramm of $\lambda$
through the main diagonal $y=-x$). The difference $\lambda_k-\lambda_{k+1}$
(where non-existing parts are considered as parts of length $0$)
equals hence the number $\nu^t_k$ of parts having length $k$
in the transposed partition $\lambda^t=(1^{\nu_1^t}\ 2^{\nu_2^t}\cdots)$
of $\lambda$. Summing over all partitions of $n$ yields then the recursion
relation $\lambda_k(n)=\nu_k(n)+\lambda_{k+1}(n)$.

The proof of the equality $\nu_k(n)=\gamma_k(n)+\nu_{k+1}(n)$ uses
generating series.
Introducing the numbers 
$$\begin{array}{l}
\displaystyle m_k(\lambda)=\sharp\{i\ \vert\ \nu_i\geq k\}\ ,\cr
\displaystyle m_k(n)=\sum_{\vert\lambda\vert=n}m_k(\lambda)\end{array}$$
one has obviously $\gamma_k(n)=m_k(n)-m_{k+1}(n)$.
We have hence to show the equality $m_k(n)=\nu_k(n)$ for $1\leq k<n$
(the equalities $m_n(n)=\nu_n(n)=1$ are easy).

We introduce the generating function 
$$\psi_k(y,z)=\sum_\lambda y^{m_k(\lambda)} z^{\vert \lambda\vert}\ .$$
One has
$$\begin{array}{l}
\displaystyle \psi_k(y,z)=\sum_{I\subset\{1,2,3,\dots\},\ \sharp(I)<\infty}
\left(\prod_{i\in I}\frac{yz^{ki}}{1-z^i}\right)\left(\prod_{1\leq j\not\in I}
\frac{1-z^{kj}}{1-z^j}\right)\cr
\displaystyle \qquad =\prod_{j=1}^\infty\left(\frac{yz^{kj}}{1-z^j}+
\frac{1-z^{kj}}{1-z^j}\right)=
\prod_{j=1}^\infty\left(\frac{1}{1-z^j}-(1-y)
\frac{z^{kj}}{1-z^j}\right)\ .
\end{array}$$
A small computation yields
$$\frac{\partial \psi_k(y,z)}{\partial y}=\psi_k(y,z)
\sum_{j=1}^\infty\frac{z^{kj}}{1-(1-y)z^{kj}}$$
and we have hence
$$\sum_{n=0}^\infty m_k(n)\ z^n=\frac{\partial \psi_k}{\partial y}(1,z)=
\left(\prod_{j=1}^\infty\frac{1}{1-z^j}\right)\ \frac{z^k}{1-z^k}=
\sum_{n=0}^\infty \nu_k(n)\ z^n$$
(cf. Theorem 1.1 for the last equality) which finishes the proof.\hfill QED

Corollary 1.4 results immediately from Theorem 1.3 and from
the last equality in Theorem 1.1.

{\bf Proof of Proposition 1.5.} A partition
$$\lambda=(\lambda_1\geq \lambda_2\geq\dots\geq \lambda_k=i\geq \lambda_{k+1}\geq\dots)$$
with $\lambda_k=i$ of $n=\sum_{j=1}\lambda_j$ can be written as
$$\lambda=(i+(\lambda_1-i)\geq i+(\lambda_2-i)\geq\dots\geq i+(\lambda_{k-1}-i)\geq i
\geq\lambda_{k+1}\geq\dots)\ .$$
Such partitions are hence in bijection with pairs of partitions
$$\begin{array}{l}
\alpha=(\alpha_1=(\lambda_1-i)\geq \alpha_2=(\lambda_2-i)\geq\dots\geq 
\alpha_{k-1}=(\lambda_{k-1}-i)
\geq 0)\ ,\cr
\omega=(\omega_1=\lambda_{k+1}\geq\omega_2=\lambda_{k+2}\geq\dots)
\end{array}$$
with $\alpha$ having at most $k-1$ non-zero parts and $\omega$ having all 
parts $\leq i$. The conjugate partition $\alpha^t$ of $\alpha$ has 
hence only parts 
$\leq k-1$. Such a pair $\alpha^t,\omega$ of partitions yields hence a unique
partition with $\lambda_k=i$ of the integer $n=ki+\sum_{j=1}^{k-1}\alpha^t_j+
\sum_{j=k+1} \omega_j$ and contributes hence with ${i\choose d}$ 
to the $k-$th
coordinate ${\lambda_k(n)\choose d}$ of ${\lambda(n)\choose d}$. 
Summing up over $i\in{\bf N}$ yields easily the generating series for
${\lambda_k(n)\choose d}$.

Considering the generating series for ${\nu_k(n)\choose d}$ one has
$$\sum_n{\nu_k(n)\choose d}z^n=\left(\sum_j{j\choose d}z^{jk}\right)
\prod_{i\not=k}\frac{1}{1-z^i}$$
and the (easy) equality
$$\sum_j{j\choose d}Z^j=\frac{1}{Z}\left(\frac{Z}{1-Z}\right)^{d+1}$$
implies the result.
\hfill QED

{\bf Proof of Theorem 1.7.} Given a partition $\lambda=(\lambda_1,\lambda_2,\dots)$
the definition 
$$\lambda^t_k=\sharp\{i\ \vert \ \lambda_i\geq k\}$$
for the $k-$th part of its transposed partition 
$\lambda^t=(\lambda^t_1,\lambda^t_2,\dots)$ shows the equalities
$$\sum_\lambda y^{\lambda_k}\ z^{\vert\lambda\vert}=
\sum_\lambda y^{(\lambda^t)_k}\ z^{\vert\lambda\vert}=
\left(\prod_{i=1}^{k-1}
\frac{1}{1-z^i}\right)\left(\prod_{j=k}^\infty
\frac{1}{1-yz^j}\right)\ .$$
Denote this series by $\varphi_k(y,z)$. An easy computation yields
$$\varphi_k(y,z)=\left(\prod_{i=1}^\infty
\frac{1}{1-z^i}\right)\prod_{j\geq k}\left(1-(y-1)\frac{z^j}{1-z^j}\right)^{-1}\ .$$
Applying the identity
$$\prod_{i}\left(1-x_i\right)^{-1}=\hbox{exp}\left(
\sum_{l=1}^\infty\frac{\sum_i x_i^l}{l}\right)$$
of formal power series to the last factor we get
$$\begin{array}{l}
\displaystyle 
\prod_{j\geq k}\left(1-(y-1)\frac{z^j}{1-z^j}\right)^{-1}
=\hbox{exp}\left(\sum_{l=1}^\infty\frac{(y-1)^l}{l}\sigma_l(k)\right)\cr
\displaystyle \quad =\sum_{n=0}^\infty \frac{(y-1)^n}{n!}\sum_{\vert(
1^{\nu_1}\cdot 2^{\nu_2}\cdots)\vert=n}\frac{n!}{\left(\sum_i \nu_i\right)!}
{\sum_i\nu_i\choose
\nu_1\ \nu_2\dots}\prod_i\left(\frac{\sigma_i(k)}{i}\right)^{\nu_i}\cr
\displaystyle \quad =\sum_{n=0}^\infty \frac{(y-1)^n}{n!}\ S_n(k)\ .\end{array}$$
We have thus
$$d!\sum_\lambda {\lambda_k\choose d}z^{\vert \lambda\vert}=
\frac{\partial^d \varphi_k}{\partial y^d}(1,z)=
\left(\prod_{i=1}^\infty \frac{1}{1-z^i}\right)
S_d(k)$$
which finishes the proof by comparison with Proposition 1.5.\hfill QED

\vskip0.5cm
We thank S. Attal and M-L. Chabanol for comments and for their 
interest in this work.

{\bf References}

{\bf [A]} G.E. Andrews, {\it The Theory of Partitions, Encyclopedia of 
Mathematics and its Applications}, Addison-Wesley (1976).
 
{\bf [EL]} P. Erd\"os, J. Lehner, {\it The distribution of the number of
summands in the partition of a positive integer}, Duke Math. Journal
{\bf 8}, 335-345 (1941).

{\bf [O]} A. Okounkov, {\it Random Matrices and Random Permutations},
IMRN {\bf 20}, 1043-1095 (2000).

{\bf [S]} Stanley, {\it Enumerative Combinatorics 2}, Cambridge 
University Press (1999).

{\bf [ST]} M. Szalay, P. Turan, {\it On some problems of the statistical
theory of partitions with application to characters of the symmetric
group I}, Acta Math. Sci. Hung. {\bf 29}, 361-379 (1977).

{\bf [VK]} A. Vershik, S. Kerov, {\it Asymptotics of the Plancherel 
measure of the symmetric group and the limit form of Young 
tableaux, Soviet Math. Dokl. {\bf 18}, 1977, 527-531.

Roland Bacher
 
INSTITUT FOURIER
 
Laboratoire de Math\'ematiques
 
UMR 5582 (UJF-CNRS)
 
BP 74
 
38402 St MARTIN D'H\`ERES Cedex (France)
 
e-mail: Roland.Bacher@ujf-grenoble.fr  

\vskip0.5cm

Laurent Manivel
 
INSTITUT FOURIER
 
Laboratoire de Math\'ematiques
 
UMR 5582 (UJF-CNRS)
 
BP 74
 
38402 St MARTIN D'H\`ERES Cedex (France)
 
e-mail: Laurent.Manivel@ujf-grenoble.fr 

\end{document}